\documentclass[12pt,leqno]{article}
\usepackage{amsfonts}
\usepackage{amsmath, amsthm, amsfonts, amssymb, color}
\usepackage{mathrsfs}
\usepackage{color}
\theoremstyle{definition}

\newcommand{\scr}[1]{\mathscr #1}
\definecolor{wco}{rgb}{0.5,0.2,0.3}

\numberwithin{equation}{section} \theoremstyle{remark}

\newcommand{\ua}{\uparrow}

\title{{\bf Killed Distribution Dependent SDE for Nonlinear Dirichlet Problem }\footnote{Supported in
 part by  NNSFC (11831014, 11921001) and the National Key R\&D Program of China (No. 2020YFA0712900).} }
\author{{\bf  Feng-Yu Wang   }\\
\footnotesize{ Center for Applied Mathematics, Tianjin University, Tianjin 300072, China}\\
  \footnotesize{  Department of Mathematics,
Swansea University, Bay Campus, SA1 8EN, United Kingdom}\\
\footnotesize{    wangfy@tju.edu.cn}}
\begin{document}
\allowdisplaybreaks
\def\R{\mathbb R}  \def\ff{\frac} \def\ss{\sqrt} \def\B{\mathbf
B}
\def\N{\mathbb N} \def\kk{\kappa} \def\m{{\bf m}}
\def\ee{\varepsilon}\def\ddd{D^*}
\def\dd{\delta} \def\DD{\Delta} \def\vv{\varepsilon} \def\rr{\rho}
\def\<{\langle} \def\>{\rangle}
  \def\nn{\nabla} \def\pp{\partial} \def\E{\mathbb E}
\def\d{\text{\rm{d}}} \def\bb{\beta} \def\aa{\alpha} \def\D{\scr D}
  \def\si{\sigma} \def\ess{\text{\rm{ess}}}\def\s{{\bf s}}
\def\beg{\begin} \def\beq{\begin{equation}}  \def\F{\scr F}
\def\Ric{\mathcal Ric} \def\Hess{\text{\rm{Hess}}}
\def\e{\text{\rm{e}}} \def\ua{\underline a} \def\OO{\Omega}  \def\oo{\omega}
 \def\tt{\tilde}\def\[{\lfloor} \def\]{\rfloor}
\def\cut{\text{\rm{cut}}} \def\P{\mathbb P} \def\ifn{I_n(f^{\bigotimes n})}
\def\C{\scr C}      \def\aaa{\mathbf{r}}     \def\r{r}
\def\gap{\text{\rm{gap}}} \def\prr{\pi_{{\bf m},\nu}}  \def\r{\mathbf r}
\def\Z{\mathbb Z} \def\vrr{\nu} \def\ll{\lambda}
\def\L{\scr L}\def\Tt{\tt} \def\TT{\tt}\def\II{\mathbb I}
\def\i{{\rm in}}\def\Sect{{\rm Sect}}  \def\H{\mathbb H}
\def\M{\mathbb M}\def\Q{\mathbb Q} \def\texto{\text{o}} \def\LL{\Lambda}
\def\Rank{{\rm Rank}} \def\B{\scr B} \def\i{{\rm i}} \def\HR{\hat{\R}^O}
\def\to{\rightarrow} \def\gg{\gamma}
\def\EE{\scr E} \def\W{\mathbb W}
\def\A{\scr A} \def\Lip{{\rm Lip}}\def\S{\mathbb S}
\def\BB{\scr B}\def\Ent{{\rm Ent}} \def\i{{\rm i}}\def\itparallel{{\it\parallel}}
\def\g{{\mathbf g}}\def\Sect{{\mathcal Sec}}\def\T{\mathcal T}\def\BB{{\bf B}}
\def\f\ell \def\g{\mathbf g}\def\BL{{\bf L}}  \def\BG{{\mathbb G}}
\def\Bd{{D^E}} \def\BdP{D^E_\phi} \def\Bdd{{\bf \dd}} \def\Bs{{\bf s}} \def\GA{\scr A}
\def\Bg{{\bf g}}  \def\Bdd{\psi_B} \def\supp{{\rm supp}}\def\div{{\rm div}}
\def\ddiv{{\rm div}}\def\osc{{\bf osc}}\def\1{{\bf 1}}\def\BD{\mathbb D}
\def\H{{\bf H}}\def\gg{\gamma} \def\n{{\mathbf n}} \def\K{{\bf K}}\def\GG{\Gamma}
\maketitle

\begin{abstract}To characterize nonlinear Dirichlet problems in an open domain, we investigate killed distribution dependent SDEs. By constructing the coupling by projection and using the Zvonkin/Girsanov transforms,    the well-posedness is proved for three different situations: 

 1) monotone case with distribution dependent noise  (possibly degenerate);

2) singular case 
with non-degenerate distribution dependent noise;  

3) singular case 
with non-degenerate distribution independent noise. \newline In the first two cases the domain is $C^2$ smooth such that the Lipschitz continuity in initial distributions is also derived, and  in the last  case the domain is arbitrary.  

  \end{abstract} \noindent
 AMS subject Classification:\  60B05, 60B10.   \\
\noindent
 Keywords:    Nonlinear Dirichlet problem, killed distribution dependent SDE,   coupling by projection, well-posedness.

 \vskip 2cm

 \section{Introduction }

The distribution dependent  stochastic differential equation  (DDSDE) is a crucial probability model characterizing the nonlinear Fokker-Planck equation. It is known as   McKean-Vlasov SDE due to \cite{MC},  and    mean field SDE for its link to mean field particle systems,  see for instance the lecture notes \cite{SN} and the survey  \cite{HRW} for the background   and recent progress on the study of DDSDEs and applications. To characterize the nonlinear Neumann problem,
the reflecting   DDSDE   has been investigated in  \cite{W21}, see also \cite{Reis} for the  convex domain case, and see the early work \cite{SN0} for the characterization on propagations of chaos.
In this paper, we consider the killed DDSDE  which in turn to describe the  nonlinear Dirichlet problem.
In this case, the distribution is restricted to the open domain and thus might be a sub-probability measure, i.e.
the total mass may be less than $1$.

Let $O\subset \R^d$ be a connected open domain with closure $\bar O$, and let  
$$\scr P_{O}:=\big\{\mu\ \text{is\ a\ measure\ on\ }O, \ \mu(O)\le 1\big\}$$ be the space of sub-probability measures on $O$ equipped with the weak topology.
Consider the following time-distribution dependent second order differential operator on $O$:
$$L_{t,\mu}:= {\rm tr}\{(\si_t\si_t^*)(\cdot,\mu)\nn^2\}+\nn_{b_t(\cdot,\mu)},\ \ t\in [0,T],\mu\in \scr P_O,$$
where $T>0$ is a fixed constant,    $\si^*$ is the transposition of $\si$, $\nn^2$ is the Hessian operator, $\nn_b:=b\cdot \nn$ is the derivative along $b$, and for some  $m\in \mathbb N$, 
$$b: [0,T]\times O\times \scr P_O\to \R^d,\ \ \si: [0,T]\times O\times \scr P_O\to \R^d\otimes\R^m$$ 
are measurable such that  
\beq\label{C} \int_0^T\d t \int_{O} \big\{|b_t(x,\mu_t)|+\|\si_t(x,\mu_t)\|^2\big\} \mu_t(\d x)<\infty,\ \ \mu=(\mu_t)_{t\in [0,T]}\in C([0,T];\scr P_O).\end{equation}

To introduce the nonlinear Dirichlet problem for $L_{t,\mu}$ on $\scr P_O$, let  $C_D^2(O)$  be the class of 
  $f\in C_b^2(\bar O)$ with Dirichlet condition $f|_{\pp O}=0,$   where $f\in C_b^2(\bar O)$ means that $f$ is a bounded $C^2$ function on $\bar O$ with bounded first and second order derivatives. 
For any   $t\in [0,T]$ and $\mu,\nu\in \scr P_O$ such that 
$$\int_O  \big\{|b_t(x,\mu)|+\|\si_t(x,\mu)\|^2\big\}\nu(\d x)<\infty,$$ define the linear functional on  $C_D^2(O)$:
$$L_{t,\mu}^{D*}\nu:\, C_D^2 (O)\ni f\mapsto (L_{t,\mu}^{D*}\nu)(f):= \int_O L_{t,\mu} f \d\nu\in \R. $$     
The corresponding nonlinear Dirichlet problem  for $L_{t,\mu}$ is  the equation 
\beq\label{FP} \pp_t\mu_t= L_{t,\mu_t}^{D*}\mu_t,\ \ t\in [0,T]\end{equation} for $\mu: [0,T]\to \scr P_O.$   
We call $\mu_\cdot \in C([0,T]; \scr P_O)$  a solution to \eqref{FP}, if
$$\mu_t(f)=\mu_0(f)+\int_0^t \mu_s(L_{s,\mu_s}f)\d s,\ \ t\in [0,T],f\in C_D^2(O),$$  
where  $\mu(f):=\int f\d\mu$ for a measure $\mu$ and $f\in L^1(\mu)$. 

When $\mu_t(\d x)=\rr_t(x)\d x$, \eqref{FP} reduces to the nonlinear Dirichlet problem
$$\pp_t \rr_t= L_{t,\rr_t}^{D*}\rr_t,\ \ t\in [0,T],$$
where $L_{t,\rr_t}:= L_{t,\rr_t(x)\d x},$ in the sense that
$$\int_O (f\rr_t)(x)\d x =\int_O (f\rr_0)(x)\d x+\int_0^t \d s\int_O  (\rr_sL_{s,\rr_s}f)(x)\d x,\ \ t\in [0,T],f\in C_D^2(O).$$
 
To characterize \eqref{FP}, we consider the following killed distribution dependent SDE on $\bar O$:
\beq\label{E1} \d X_t =\1_{\{t<\tau(X)\}}\big\{b_t(X_t, \L_{X_t}^O)\d t +\si_t(X_t,\L_{X_t}^O)\d W_t\big\},\ \ t\in [0,T],\end{equation}
where $\1$ is the indicated function, $W_t$ is the $m$-dimensional Brownian motion on a complete filtration probability space $(\OO,\{\F_t\}_{t\ge 0},\P)$,
$$\tau(X):=  \inf\{t\in [0,T]:  X_t\in\pp O\}$$
with   $\inf\emptyset =\infty$ by convention,
and for an $\bar O$-valued random variable $\xi$,
$$\scr L^O_\xi:=\P(\xi\in O\cap\cdot)$$
is the distribution of $\xi$ restricted to $O$, which we call the $O$-distribution of $\xi$. When different probability spaces are concerned, we denote $\L_\xi^O$ by $\L_{\xi|\P}^O$  to emphasize the reference probability measure. 

\beg{defn}\label{D1.1} A continuous adapted process $(X_t)_{t\in [0,T]}$ on $\bar O$ is called a solution of \eqref{E1}, if $\P$-a.s.
$$\int_0^{T\land \tau(X)} \big\{|b_t(X_t,\L^O_{X_t})|+\|\si_t(X_t,\L_{X_t}^O)\|^2\big\}\d t<\infty$$ and
$$X_t=X_0+\int_0^{t\land \tau(X)} \big\{b_s(X_s, \L^O_{X_s})\d s+ \si_s(X_s,\L^O_{X_s})\d W_s\big\},\ \ t\in [0,T].$$
We call $(\tt X_t, \tt W_t)$ a weak solution to \eqref{E1}, if there exists a complete filtration probability space
 $(\tt\OO,\{\tt\F_t\}_{t\in [0,T]},\tt\P)$ such that $\tt W_t$ is $m$-dimensional Brownian motion and $\tt X_t$ solves \eqref{E1} for $\tt W_t$ replacing $W_t$.
\end{defn}

\paragraph{Remark 1.1.}    (1) It is easy to see that for any (weak) solution $X_t$ of \eqref{E1}, $\mu_t:=\L_{X_t}^O$ solves the nonlinear Dirichlet problem \eqref{FP}. Indeed, since $\d X_t=0$ for $t\ge \tau(X)$, we have
$$X_t=X_{\tau(X)}\in \pp O,\ \ t\ge \tau(X),$$
so that 
$$X_t= X_{t\land\tau(X)},\ \ \L^O_{X_t}(\d x)=\P(t<\tau(X), X_t\in\d x),\ \ t\in [0,T].$$
By this and It\^o's formula, for any $f\in C_D^2(O)$ we have 
\beg{align*} &  \mu_t(f)  =\E[(\1_{O}f)(X_t)]= \E [f(X_t)]\\
& = \E[f(X_0)]+\E \int_0^{t}\1_{\{s<\tau(X)\}} L_{s,\mu_s}f(X_s)\d s\\
&= \mu_0(f) +\int_0^t\mu_s( L_{s,\mu_s} f)\d s,   \ \ t\in [0,T].\end{align*}

(2) An alternative model to \eqref{E1} is 
\beq\label{E1'}  \d X_t =\1_{O}(X_t)  \big\{b_t(X_t, \L_{X_t}^O)\d t +\si_t(X_t,\L_{X_t}^O)\d W_t\big\},\ \ t\in [0,T].\end{equation} 
A solution of \eqref{E1} also solves  \eqref{E1'};  while  for a    solution $X_t$ to \eqref{E1'},
$$\tt X_t:= X_{t\land \tau(X)}$$ solves \eqref{E1}. In general, a solution of \eqref{E1'} does not have to solve \eqref{E1}. For instance, let $d=m=1$ and $O=(0,\infty),$
 consider $\si_t(x,\mu)=2x, b_t(x,\mu)=2\ss x.$ Let $Y_t$ solve  the   SDE
 $$\d Y_t= Y_t\d W_t +\Big(1-\ff 1 2 Y_t\Big)\d t,\ \ Y_0=0.$$ Then    $X_t:=(Y_t)^2$   solves \eqref{E1'} but does not solve \eqref{E1},  since $\tau(X)=0$ and  $X_t>0$ (i.e. $X_t\notin\pp O$) for $t> 0$.  See   \cite{ZZ} for the study of \eqref{E1'} for $\si_t(x,\mu)=\si_t(x)$ independent of $\mu$.

(3)  The SDE \eqref{E1'} can be formulated as the usual DDSDE on $\R^d$, so that the superposition principle in \cite{BR2} applies. 
 More precisely, let  $\scr P$ be the space of probability measures on $\R^d$, and define 
$$\bar b_t(x,\mu):= \1_O(x) b_t(x, \mu(O\cap\cdot)),\ \ \bar \si_t(x,\mu):= \1_O(x) \si_t(x, \mu(O\cap\cdot))$$
for $(t,x,\mu)\in [0,T]\times \R^d\times \scr P.$  Then \eqref{E1'} becomes the following  DDSDE on $\R^d$:
$$  \d X_t =\bar b_t(X_t, \L_{X_t})\d t +\bar \si_t(X_t,\L_{X_t})\d W_t,\ \ t\in [0,T].$$

 \
 
We often  solve  \eqref{E1} for $O$-distributions  in  a non-empty sub-space $\hat{\scr P}_O$  of $\scr P_O$, which is equipped with the weak topology as well.  

\beg{defn} (1) If for any $\F_0$-measurable random variable $X_0$ on $\bar O$ with $\L_{X_0}^O\in \hat{\scr P}_O$,  \eqref{E1} has a unique solution starting at $X_0$ such that 
$\L_{X}^O:= (\L_{X_t}^O)_{t\in [0,T]}\in C([0,T];\hat{\scr P}_O)$, we call the SDE   strongly well-posed for $O$-distributions in $\hat{\scr P}_O$. 

 (2)   We call the SDE weakly unique for    $O$-distributions in $\hat{\scr P}_O$,   if for any two weak solutions $(X_t^i,W_t^i)$ w.r.t. $(\OO^i, \{\F_t^i\}_{t\in [0,T]},\P^i) (i=1,2)$ with 
 $\L_{X_0^1|\P^1}^O=\L_{X_0^2|\P^2}^O\in\hat{\scr P}_O$, we have  $\L_{X^1|\P^1}^O= \L_{X^2|\P^2}^O.$ 
 We call \eqref{E1} weakly well-posed for $O$-distributions in $\hat{\scr P}_O$, if for any initial $O$-distribution $\mu_0\in \hat{\scr P}_O$, it has a unique weak solution for    $O$-distributions in $\hat{\scr P}_O.$
 
 (3) The SDE  \eqref{E1} is called well-posed for $O$-distributions in $\hat{\scr P}_O$, if it is both strongly and weakly well-posed for $O$-distributions in $\hat{\scr P}_O$.    \end{defn} 
  
  When \eqref{E1} is well-posed for $O$-distributions in $\hat{\scr P}_O,$  for any $\mu\in \hat{\scr P}_O$ and $t\in [0,T]$, let 
  $$P_t^{D*}\mu=\L_{X_t}^O,\ \ t\in [0,T], \ \L_{X_0}^O=\mu.$$    
  
We will study the well-posedness under the following assumption.

 \beg{enumerate}\item[{\bf (H)}]    For any  $\mu\in C([0,T];\hat{\scr P}_O)$,  the killed SDE
\beq\label{EM}  \d X_t^\mu = \1_{\{t<\tau(X^\mu)\}}\big\{b_t(X_t^\mu, \mu_t)\d t +\si_t(X_t^\mu,\mu_t)\d W_t\big\},\ \ t\in [0,T]\end{equation}
 is well-posed for initial value $X_0^\mu$ with $\L_{X_0^\mu}^O=\mu_0$, and $\L_{X^\mu}^O\in C([0,T];\hat{\scr P}_O).$
\end{enumerate}
 Under this assumption, we define a map
\beq\label{SB0} C([0,T];\hat{\scr P}_O)\ni\mu\mapsto \Phi\mu:= \L_{X^\mu}^O:=( \L_{X_t^\mu}^O)_{t\in [0,T]}\in C([0,T];\hat{\scr P}_O).\end{equation}
 It is clear that a solution of \eqref{EM} solves \eqref{E1} if and only if $\mu$ is a fixed point of $\Phi$. So, we have the following result.

\beg{thm} \label{T1} Assume {\bf (H)}. If for any $\gg\in \hat{\scr P}_O$, $\Phi$ has a unique fixed point in 
$$\C^\gg:=\{\mu\in   C([0,T];\hat{\scr P}_O), \mu_0=\gg\},$$ then   $\eqref{E1}$ is well-posed for  $O$-distributions in $\hat{\scr P}_O. $
\end{thm}

In the remainder of the paper, we apply Theorem \ref{T1}  to  following three different situations: 
\beg{enumerate} \item The monotone case with distribution dependent noise (possibly degenerate);
\item The singular case with non-degenerate distribution dependent noise;
\item The singular case with non-degenerate distribution independent noise.\end{enumerate}
In the first two situations,  we need $O$ to be $C^2$-smooth to apply the coupling by projection, and the coefficients are local Lipschitz continuous in distributions with respect to the $L^1$ or truncated $L^1$ Wasserstein distance. In the last  case, the domain is arbitrary,  and  $b_t(x,\cdot)$ is only  local Lipschitz continuous in   a weighted variation distance, but the noise has to be distribution independent, i.e. $\si_t(x,\mu)=\si_t(x)$.

\section{Monotone case  } 

In this part, we solve \eqref{E1}   under   monotone conditions  with respect to  the $L^1$ or     truncated $L^1$ Wasserstein distances:
\beq\label{W1}\beg{split}&  \W_1(\mu,\nu):=\inf_{\pi\in \C_O(\mu,\nu)} \int_{O\times O}   |x-y|\pi(\d x,\d y), \\
& \hat \W_1(\mu,\nu):=\inf_{\pi\in \C_O(\mu,\nu)} \int_{O\times O} (1\land |x-y|)\pi(\d x,\d y),\ \ \mu,\nu\in \scr P_O,\end{split}\end{equation} 
where $\pi\in \C_O(\mu,\nu)$ means that $\pi$ is a probability measure on $\bar O\times \bar O$ such that 
$$\pi(\{\cdot\cap O\}\times \bar O)=\mu,\ \ \pi(\bar O\times \{\cdot\cap O\})=\nu.$$
  
\subsection{Monotonicity in $\hat \W_1$}

 \beg{enumerate}  \item[$(A_1)$]    
 For any $\mu\in C([0,T];\scr P_O)$, $b_t(x,\mu_t)$ and $\si_t(x,\mu_t)$ are continuous in $x\in O$ such that for any $N\ge 1$ and $O_N:= \{x\in O: |x|\le N\}$,  
 $$\int_0^T \sup_{O_N} \big\{|b_t(\cdot,\mu_t)|+\|\si_t(\cdot,\mu_t)\|^2\big\}\d t <\infty.$$
Moreover, there exists $K\in L^1([0,T]; (0,\infty))$  such that for any   $x,y\in O$ and $\mu,\nu\in \scr P_O$,  
\beg{align*}  &2\<b_t(x,\mu)-b_t(y,\nu),x-y\>+\|\si_t(x,\mu)-\si_t(y,\nu)\|_{HS}^2\le K(t)\big\{|x-y|^2+ \hat \W_1(\mu,\nu)^2\big\},\\
&2\<b_t(x,\mu),x\>+\|\si_t(x,\mu)\|_{HS}^2\le K(t) \big(1+|x|^2\big),\ \  \ t\in [0,T].\end{align*}   \end{enumerate} 
 \beg{enumerate} 
\item[$(A_2)$]  There exists $r_0\in (0,1]$ such that the distance function $\rr_\pp$ to $\pp O$ is $C^2$-smooth in 
$$\pp_{r_0}O:=\{x\in \bar O:\ \rr_\pp(x)\le r_0\big\},$$
and there exists a constant  $\aa>0$    such that 
$$ |\si_t(x,\mu)^*\nn\rr_\pp(x)|^{-2}\le  \aa,\ \ \ L_{t,\mu} \rr_\pp(x)\le  \aa,\ \ \ x\in \pp_{r_0}O, t\in [0,T].$$ \end{enumerate}

 \beg{thm}\label{T2.1} Assume   $(A_1)$ and  $(A_2)$. Then the following assertions hold. 
 \beg{enumerate} \item[$(1)$] $\eqref{E1}$ is well-posed for $O$-distributions in $\scr P_O$.
Moreover,  for any $p\ge 1$  there exists a constant $c  >0$ such that for any solution  $X_t$  to $\eqref{E1}$  for $O$-distributions in $ \scr P_O$,  
 \beq\label{A01} \E \Big[\sup_{t\in [0,T]} |X_t|^p \Big|\F_0\Big]\le c \big(1+|X_0|^p\big).\end{equation} 
 \item[$(2)$] There exists a constant $c>0$ such that   
\beq\label{A02} \sup_{t\in [0,T]} \hat\W_1(P_t^{D*}\mu, P_t^{D*}\nu) \le c  \hat\W_1(\mu,\nu),\ \ \mu,\nu\in \scr P_O.\end{equation}\end{enumerate} \end{thm}

Under assumption    $(A_1)$,  for any $\mu\in C([0,T];\scr P_O)$   the SDE  \eqref{EM}  satisfies the semi-Lipschitz condition before the hitting time 
$\tau(X^\mu)$,   hence it is well-posed   
  and  for any $p\ge 1$ there exists a constant $c>0$ uniformly in $\mu$ such that   
 \beq\label{MES}   \E\Big[\sup_{t\in [0,T]}|X_{t\land\tt\tau}^\mu|^p\Big|\F_0\Big] =  \E\Big[\sup_{t\in [0,T]}|X_{t\land\tau(X^\mu)\land\tt\tau}^\mu|^p\Big|\F_0\Big]  \le c (1+|X_0^\mu|^p)\end{equation} holds for any solution $X_t^\mu$ of \eqref{EM} and any stopping time $\tt\tau$. 
  
  By Theorem \ref{T1}, to prove the well-posedness of 
 \eqref{E1} for $O$-distributions in $\scr P_O$,  it remains to show  that  for any $\gg\in \scr P_O$, the map 
 $$\Phi\mu:= \L_{X^\mu}^O=(\L_{X_t^\mu}^O)_{t\in [0,T]},\ \ \mu\in C([0,T];\scr P_O)$$ 
 has a unique fixed point in 
 $$\C^\gg:= \{\mu\in C([0,T]; \scr P_O): \mu_0=\gg\}.$$    To this end, for $i=1,2,$ 
let  $\mu^i\in C([0,T];\scr P_O)$, and  let $X_t^i$ solve \eqref{EM} for $\mu^i$ replacing $\mu$ with $\L_{X_0^i}^O=\mu_0^i$, i.e.
\beq\label{Mi}   \d X_t^i = \1_{\{t<\tau(X^i)\}}\big\{b_t(X_t^i, \mu_t^i)\d t +\si_t(X_t^i,\mu_t^i)\d W_t\big\},\ \ t\in [0,T], \L_{X_0^i}^O=\mu_0^i.\end{equation}
Simply denote
$$\tau_i= \tau(X^i)\ \text{for}\ i=1,2,\ \ \tau_{1,2}:=\tau_1\land\tau_2.$$ 

Since 
$$\GG:=\big\{(x,y):\ x\in O, y\in \pp O, |x-y|=\rr_\pp(x)\big\}$$
is a measurable subset of $O\times\pp O$ and $\GG_x:=\{y\in\pp O: (x,y)\in\GG\}\ne \emptyset$ for any $x\in O$,  by the measurable selection theorem (see \cite[Theorem 1]{Ev}),
there exists a measurable map  $P_\pp: O\to\pp O$   such that  
\beq\label{PR} |P_\pp x-x|=\rr_\pp(x),\ \ x\in O.\end{equation}    
  We  will use  the following coupling by projection.   

\beg{defn}  \label{D2.1} 
 The   \emph{coupling by projection} $(\bar X_t^1,\bar X_t^2)$ for  $(X_t^1,X_t^2)=(X_{t\land\tau_1}^1,X_{t\land\tau_2}^2)$ is defined as 
 \beq\label{CP} (\bar X_t^1,\bar X_t^2):=\beg{cases} (X_t^1,X_t^2),\ &\text{if} \ t\le \tau_{1,2},\\
 (X_t^1, P_\pp X_t^1),\ &\text{if} \ \tau_2< t\land \tau_1, \\
  (P_\pp X_t^2,   X_t^2), &\text{otherwise.}  \end{cases}\end{equation} 
 \end{defn} 
 
  It is easy to see that $\L_{\bar X_t^i}^O= \L_{X_t^i}^O=(\Phi \mu^i)_t$ for $i=1,2;$ i.e. the distribution $\L_{(\bar X_t^1,\bar X_t^2)}$ of the coupling by projection $(\bar X_t^1,\bar X_t^2)$ satisfies 
 $$\L_{(\bar X_t^1,\bar X_t^2)}\in \C_O((\Phi \mu^1)_t,( \Phi \mu^2)_t).$$ Thus, by \eqref{W1} and Definition \ref{D2.1}, 
\beq\label{PW} \beg{split} &\hat \W_1((\Phi\mu^1)_t, (\Phi\mu^2)_t)\le \E\big[1\land |\bar X_t^1-\bar X_t^2|\big]
\le \E\big[1\land |X_{t\land\tau_{1,2}}^1-\bar X_{t\land\tau_{1,2}}^2|\big]\\
&\qquad +r_0^{-1}\E[\{r_0\land \rr_\pp(X_t^1)\}\1_{\{t\land\tau_1\ge\tau_2\}}\big] +r_0^{-1} \E[\{r_0\land \rr_\pp(X_t^2)\}\1_{\{t\land\tau_2\ge\tau_1\}}\big].\end{split} \end{equation}

 \beg{lem}\label{L2.1} Assume $(A_1)$.  
Then there exists a constant $c>1$ such that for any $t\in [0,T]$ and $\mu^1,\mu^2\in C([0,T];\scr P_O),$ 
\beq\label{C1} \E\big[  |X_{t\land\tau_{1,2}}^1-X_{t\land\tau_{1,2}}^2|^2\big|\F_0\big]\le c\,   |X_0^1-X_0^2|^2+ c \int_0^t K(s)  \hat \W_1(\mu_s^1,\mu_s^2)^2 \d s.\end{equation}
Consequently, for any $t\in [0,T]$, 
\beq\label{C2}  \E\big[ 1\land  |X_{t\land\tau_{1,2}}^1-X_{t\land\tau_{1,2}}^2|\big]\le \ss c\, \E[1\land |X_0^1-X_0^2|] +\bigg( c \int_0^t K(s)  \hat \W_1(\mu_s^1,\mu_s^2)^2 \d s\bigg)^{\ff 1 2}.\end{equation}
 \end{lem}
\beg{proof} It suffices to prove \eqref{C1}, which implies \eqref{C2} due to  Jensen's inequality.  

By $(A_1)$ and It\^o's formula, we obtain 
$$\d |X_t^1-X_t^2|^2\le K(t) \big\{|X_t^1-X_t^2|^2+\hat \W_1(\mu_t^1,\mu_t^2)^2\big\}\d t + \d M_t,\ \ t\in [0,T\land \tau_{1,2}]$$ 
for some local martingale $M_t$. This and \eqref{MES} imply that   
$$\bb_t:=  \E\big[  |X_{t\land\tau_{1,2}}^1-X_{t\land\tau_{1,2}}^2|^2\big|\F_0\big] $$ is bounded in $t\in [0,T] $ and satisfies
$$\bb_t\le \bb_0 + \int_0^t K(s)\big\{\bb_s +\hat \W_1(\mu_s^1,\mu_s^2)^2\big\}\d s,\ \ t\in [0,T].$$
By Gronwall's inequality, we prove \eqref{C1}. 
\end{proof}

 \beg{lem}\label{L2.2} Assume $(A_2)$. Then  there exists a constant $c>1$   independent of $\mu$  such that  for any solution $X_t^\mu$ to $\eqref{EM}$ and any stopping time $\tt\tau$,
$$ \1_{\{t\land\tau(X^\mu)\ge \tt\tau\}}\E \big[r_0\land \rr_\pp(X_{t}^\mu)\big| \F_{\tt\tau} \big]\le c   \1_{\{t\land\tau(X^\mu)\ge \tt\tau\}}\rr_\pp(X_{t\land\tt\tau}^\mu),\ \ t\in [0,T].$$
 
 \end{lem}
 
 \beg{proof}  By the strong Markov property of $X_t^\mu$ which is implied by the well-posedness of \eqref{EM}, we may and do assume that $\tt\tau=0$ and   $x=X_0^\mu\in O$, such that the desired estimate becomes
 \beq\label{*0} \E^x \big[r_0\land \rr_\pp(X_{t}^\mu) \big]\le c  \rr_\pp(x),\ \ t\in [0,T],\end{equation} 
 where $\E^x$ is the expectation under the  probability $\P^x$ for $X_t^\mu$ starting at $x$.  
 If $\rr_\pp(x)\ge \ff {r_0}4$, this inequality holds for $c:= 4$. So, it suffices to prove for $\rr_\pp(x)<\ff {r_0}4.$ 
 
 Let $h\in C^\infty([0,\infty))$ such that 
 $$h'\ge 0,\ \ h''\le 0,\ h(r)=r\ \text{for}\ r\in [0,r_0/2],\ \ h'(r)=0\ \text{for}\ r\ge r_0.$$  
 By $(A_2)$, 
 \beq\label{*1} \d h(\rr_\pp(X_t^\mu)) \le   \aa  \d t+ \d M_t,\ \ t\in [0,T\land \tau(X^\mu)],\end{equation} 
where $M_t$ is a martingale with
\beq\label{*N} \d\<M\>_t\ge \aa^{-1}  \d t,\ \ t\le \hat \tau:=\inf\{t\ge 0: \rr_\pp(X_t^\mu)\ge r_0/2\}.\end{equation} 
By \eqref{*1} we obtain 
 \beq\label{*2}  \E^x[r_0\land \rr_\pp(X_t^\mu)]\le  2 \E^x[h(\rr_\pp(X_{t\land\tau(X^\mu)}^\mu))]\le 2 \rr_\pp(x)+ 2\aa\E^x[t\land\tau(X^\mu)].\end{equation} 
 On the other hand, let 
 $$\eta_t:= \int_0^{\rr_\pp(X_t^\mu)} \e^{-2\aa^2 s}\d s \int_s^{r_0} \e^{2\aa^2\theta}\d \theta,\ \ 
t\in [0,T\land \tau(X^\mu)\land\hat\tau].$$
Since $h(r)=r$ for $r\le \ff{r_0} 2,$ by \eqref{*1}, \eqref{*N} and It\^o's formula, we find a martingale $\tt M_t$ such that  
$$\d \eta_t\le - \d t +\d\tt M_t,\ \ t\in [0,T\land \tau(X^\mu)\land\hat\tau].$$
  Consequently,
\beq\label{LOT1} \E^x   [t\land \tau(X^\mu)\land\hat\tau] \le     \eta_0  \le c_1 \rr_\pp(x)\end{equation}
 holds for some constant $c_1>0$.   Therefore,
\beq\label{B*}\beg{split}& \E^x  [t\land \tau(X^\mu)]\le \E^x[t\land\tau(X^\mu)\land\hat\tau]+ T\,\E^x\big[\1_{\{t\land\tau(X^\mu)>\hat\tau\}}  \big] \\
&\le c_1 \rr_\pp(x)+T\, \P^x\big(t\land\tau(X^\mu)>\hat\tau\big), \ \ t\in [0,T].\end{split}\end{equation} 
 To estimate the second term, let
 $$\xi_t:=\int_0^{\rr_\pp(X_t^\mu)} \e^{-2\aa^2 s  }\d s ,\ \  t\in [0,T\land \tau(X^\mu)\land\hat\tau].$$
By $h(r)=r$ for $r\in [0,\ff{r_0}2]$, \eqref{*1}, \eqref{*N} and It\^o's fomrula, we see that $\xi_t$ is a sup-martingale, so that 
\beq\label{B*2} \rr_\pp(x) \ge  \xi_0 \ge \E^x[\xi_{t\land\tau(X^\mu)\land\hat\tau}] \ge \P^x\big(t\land\tau(X^\mu)\ge \hat\tau\  \big)  \int_0^{r_0/2} \e^{-2\aa^2s}\d s.\end{equation}
 Combining this with \eqref{*2} and \eqref{B*}, we prove \eqref{*0} for some constant $c>0$. 
   \end{proof} 
 
 \beg{proof}[Proof of Theorem \ref{T2.1}] (a) Well-posedness.  Let $\gg:=\L_{X_0}^O$, and consider 
\beq\label{CG} \C^\gg:=\big\{\mu\in C([0,T];\scr P_O):\ \mu_0=\gg\big\}.\end{equation}
We intend to prove that $\Phi$ is contractive  in $\C^\gg$ under the  complete metric
 $$\hat\W_{1,\theta}(\mu^1,\mu^2):=\sup_{t\in [0,T]}\e^{-\theta t} \hat\W_1(\mu_t^1,\mu_t^2)$$
 for large enough $\theta>0$. Then $\Phi$ has a unique fixed point in $\C^\gg$, so that the well-posedness follows from Theorem \ref{T1}. 
 
 To this end, let $\mu^i\in \C^\gg$ and let $X_t^i$ solve \eqref{EM} with $\mu=\mu^i$ and $X_0^i=X_0, i=1,2.$  By $r_0\le 1$,  
  Lemma \ref{L2.2}, and noting that 
   $$\1_{\{t\land\tau_2\ge \tau_1\}}\rr_\pp(X_{t\land\tau_{1,2}}^2)\le \1_{\{t\land\tau_2\ge \tau_1\}} |X_{t\land\tau_{1,2}}^2-X_{t\land\tau_{1,2}}^1|,$$  we obtain
\beq\label{ORR} \beg{split}  &\E\Big[\1_{\{t\land\tau_2\ge \tau_1\}}\big\{r_0\land  \rr_\pp (X_{t\land \tau_2}^2) \big\} \Big]= \E\Big(\1_{\{t\land\tau_2\ge \tau_1\}}\E\big[\big\{r_0\land  \rr(X^2_{t\land \tau_2} )\big\}\big|\F_{\tau_1}\big] \Big)\\
 &\le c\ \E\Big[\1_{\{t\land\tau_2\ge \tau_1\}}\big\{r_0\land \rr_\pp(X_{t\land\tau_{1,2}}^2)\big\}\Big]
 \le c\  \E\big[1 \land |X_{t\land\tau_{1,2}}^1-X_{t\land\tau_{1,2}}^2|\big].\end{split}\end{equation}
 By symmetry, the same estimate holds for $\E\Big[\1_{\{t\land\tau_1\ge \tau_2\}}\big\{r_0\land  \rr_\pp (X_{t\land \tau_2}^1) \big\} \Big]$. Combining these with  $X_0^1=X_0^2=X_0$,  \eqref{PW} and \eqref{C2}, we find a constant $c_1>0$ such that  
$$\hat \W_1((\Phi\mu^1)_t, (\Phi\mu^2)_t)\le c_1 \bigg(\int_0^tK(s) \hat \W_1(\mu_s^1,\mu_s^2)^{2}\d s\bigg)^{\ff 1 2},\ \ t\in [0,T].$$
This implies that $\Phi$ is contractive in $\hat\W_{1,\theta}$ for large enough $\theta>0$.

(b) Estimate \eqref{A01}. Let $\mu_t=\L_{X_t}^O$ for the unique solution of \eqref{E1}, we have $X_t=X_t^\mu$ since $\mu$ is a fixed point of $\Phi$. So,  
\eqref{A01} follows from  \eqref{MES}. 

(c) Estimate  \eqref{A02}.  Take $X_0^1,X_0^2$ such that 
\beq\label{X0} \L_{X_0^1}^O= \mu,\ \ \L_{X_0^2}^O=\nu,\ \ \E[1\land |X_0^1-X_0^2|]=\hat\W_1(\mu,\nu).\end{equation}  Let $X_t^1$ and $X_t^2$ solve \eqref{E1}. Then they solve $\eqref{Mi}$ with
$$\mu_t^1:=\L_{X_t^1}^O= P_t^{D*}\mu,\ \ \mu_t^2:=\L_{X_t^2}^O=P_t^{D*}\nu,$$ so that  $\mu_t^i= (\Phi\mu^i)_t,\ t\in [0,T], \ i=1,2.$
Thus, by 
  \eqref{PW}, \eqref{C1} and Lemma \ref{L2.2}, we find a constant $c_2>0$ such that 
\beg{align*} &\hat \W_1(P_t^{D*}\mu, P_t^{D*}\nu) = \hat \W_1((\Phi\mu^1)_t, (\Phi\mu^2)_t)  \\
&\le c_2 \hat \W_1(\mu,\nu) + \bigg(c_2 \int_0^tK(s) \hat \W_1(P_s^{D*}\mu, P_s^{D*}\nu)^2\d s\bigg)^{\ff 1 2},\ \ t\in [0,T].\end{align*}
By Gronwall's inequality, we prove \eqref{A02} for some constant $c>0.$
 \end{proof} 
 
 \subsection{Monotonicity in $\W_1$} 
 Let $\scr P_O^1=\{\mu\in \scr P_O, \|\mu\|_1:=\mu(|\cdot|)<\infty\}.$ Define
 $$\|\mu\|_{1,T}:=\sup_{t\in [0,T]}\|\mu_t\|_1,\ \ \ \mu\in C([0,T]; \scr P_O^1).$$
 
 \beg{enumerate} 
\item[$(B_1)$] For any $\mu\in C([0,T];\scr P_O^1)$, $b_t(x,\mu_t)$ and $\si_t(x,\mu_t)$ are continuous in $x\in O$ such that for any  $N\ge 1$ and $O_N:= \{x\in O: |x|\le N\}$, 
$$\int_0^T \sup_{O_N} \big\{|b_t(\cdot,\mu_t)|+\|\si_t(\cdot,\mu_t)\|^2\big\}\d t <\infty.$$
Moreover, there exists $K\in L^1([0,T]; (0,\infty))$ such that for any   $x,y\in O$ and $\mu,\nu\in \scr P_O^1$,  
\beg{align*}  &2\<b_t(x,\mu)-b_t(y,\nu),x-y\>+\|\si_t(x,\mu)-\si_t(y,\nu)\|_{HS}^2\le K(t)\big\{|x-y|^2+ \W_1(\mu,\nu)^2\big\},\\
&2\<b_t(x,\mu),x\>+\|\si_t(x,\mu)\|_{HS}^2\le K(t) \big\{1+|x|^2+ \|\mu\|_1^2\big\},\ \  \ t\in [0,T].\end{align*} 
\item[$(B_2)$]  There exists $r_0>0$ such that  $\rr_\pp \in C^2(\pp_{r_0}O)$, and there exists an  increasing function $\aa: [0,\infty)\to [1,\infty)$    such that 
\beq\label{LU}   |\si_t(x,\mu)^*\nn\rr_\pp|^{-2}\le \aa(\|\mu\|_1),\ \ \ L_{t,\mu} \rr_\pp(x)\le   \aa(\|\mu\|_1),\ \ \ x\in \pp_{r_0}O,\end{equation} 
\beq\label{LU2} \beg{split}&2\<b_t(x,\mu),x-y\> +\|\si_t(x,\mu)\|_{HS}^2\\
 &\le K(t)\aa(\|\mu\|_1) \big(1+|x-y|^2\big),\ \  \ t\in [0,T], y\in \pp O, x\in O.\end{split}\end{equation} 
\end{enumerate} 

 \beg{thm}\label{T2.2} Assume $(B_1)$ and $(B_2)$. Then the following assertions hold. 
 \beg{enumerate} \item[$(1)$] $\eqref{E1}$ is well-posed for $O$-distributions in $\scr P_O^1$.
Moreover,  for any $p\ge 1$  there exists a constant $c>0$ such that for any solution  $X_t$  to $\eqref{E1}$ for $O$-distributions in $\scr P_O^1$,  
 \beq\label{A01'}  \E \Big[\sup_{t\in [0,T]}  |X_t|^p \Big|\F_0\Big]\le c\big(1+|X_0|+\E[\1_{O}(X_0)|X_0|] \big)^p.\end{equation} 
 \item[$(2)$] If $\aa$ is bounded, 
 then there exists a constant $c>0$ such that 
\beq\label{A02'} \sup_{t\in [0,T]}\W_1(P_t^{D*}\mu, P_t^{D*}\nu) \le  c \,  \W_1(\mu,\nu),\ \ \mu,\nu\in \scr P_O^1.\end{equation}\end{enumerate} \end{thm}

  It is standard that $(B_1)$ and $(B_2)$ imply the well-posedness of  \eqref{EM} for    $\mu\in C([0,T];\scr P_O^1)$,  and instead of \eqref{MES}, for any $p\ge 1$ there exists a constant 
  $c>0$ such that    
 \beq\label{MES'} \E\Big[\sup_{t\in [0,T]} |X_{t}^\mu|^p \Big|\F_0\Big]\le c \big(1+|X_0^\mu|^p\big)+c \int_0^t K(s) \|\mu_s\|_{1}^p\d s,\ \ t\in [0,T],\mu\in C([0,T];\scr P_O^1).\end{equation}  
 Let  $\mu^i\in C([0,T];\scr P_O^1), i=1,2$,    let $X_t^i$ solve \eqref{EM} for $\mu^i$ replacing $\mu$ with $\L_{X_0^i}^O=\mu_0^i$, and denote as before 
$$\tau_i:= \tau(X^i)\ \text{for}\ i=1,2,\ \ \tau_{1,2}:=\tau_1\land\tau_2.$$ 
Using $(B_1)$ replacing $(A_1)$, the proof of \eqref{C1} leads to
\beq\label{C1'} \E\big[  |X_{t\land\tau_{1,2}}^1-X_{t\land\tau_{1,2}}^2|^2\big|\F_0\big]\le c   |X_0^1-X_0^2|^2+ c \int_0^t K(s)   \W_1(\mu_s^1,\mu_s^2)^2 \d s,\ \ t\in [0,T],\end{equation}
and instead of \eqref{PW}, we have
 \beq\label{PW'} \beg{split} & \W_1((\Phi\mu^1)_t, (\Phi\mu^2)_t)\le \E\big[  |\bar X_t^1-\bar X_t^2|\big]
\le \E\big[  |X_{t\land\tau_{1,2}}^1-\bar X_{t\land\tau_{1,2}}^2|\big]\\
&\qquad + \E[  \rr_\pp(X_t^1)\1_{\{t\land\tau_1\ge \tau_2\}}\big] +  \E[  \rr_\pp(X_t^2)\1_{\{t\land\tau_2\ge \tau_1\}}\big].\end{split} \end{equation}
The following lemma is analogous to Lemma \ref{L2.2}.

 \beg{lem}\label{L2.2'} Assume $(B_2)$. Then  there exists an increasing function $\psi: [0,\infty)\to (0,\infty)$ which is bounded if so is $\aa$,     such that  for any $\mu\in C([0,T];\scr P_O^1)$ and any solution $X_t^\mu$ to $\eqref{EM}$ and any stopping time $\tt\tau$,
$$ \1_{\{t\land\tau(X^\mu)\ge \tt\tau\}}\E \big[  \rr_\pp(X_{t}^\mu)\big| \F_{\tt\tau} \big]\le \1_{\{t\land\tau(X^\mu)\ge \tt\tau\}}\psi(\|\mu\|_{1,T})  \rr_\pp(X^\mu_{\tt\tau}).$$
  \end{lem}

 \beg{proof} By the strong Markov property, we may assume that $\tt\tau=0$ and $x=X_0^\mu\in O$, so that it suffices to prove 
 \beq\label{*1N} \GG_t(x):=\E^x[\rr_\pp (X_t^\mu)]\le \psi(\|\mu\|_{1,T})   \rr_\pp(x),\ \ x\in O, t\in [0,T].\end{equation} 
(a) Let $\rr_\pp(x)\ge \ff{r_0}2$ and $y\in\pp O$ such that $\rr_\pp(x)=|y-x|. $
 By \eqref{LU2},   we have 
 $$\d |X_t^\mu-y|^2\le K(t)\aa(\|\mu\|_{1,T}) \big(1+|X_t^\mu-y|^2 \big)\d t +\d M_t,\ \ t\in [0,T\land\tau(X^\mu)]$$
   for some martingale $M_t$. Combining this with  $|x-y|=\rr_\pp(x),$ we obtain 
 \beg{align*}& \E^x[|X_t^\mu-y|^2]\le  \rr_\pp(x)^2+ \aa(\|\mu\|_1) \int_0^t K(s)\d s  \\
 &+   \int_0^t K(s) \aa(\|\mu\|_{1,T})\E^x[|X_s^\mu-y|^2]\d s,\ \ t\in [0,T].\end{align*}
   By Gronwall's inequality and $\rr_\pp(x)\ge \ff{r_0}2$, we find an increasing function $\psi_1: [0,\infty)\to (0,\infty)$  which is bounded if so is $\aa$, such that  
\beg{align*} \E^x[|X_t^\mu-y|^2]&\le \bigg\{ \rr_\pp(x)^2+ \aa(\|\mu\|_{1,T})\int_0^T K(s)\d s\bigg\}\e^{  \aa(\|\mu\|_{1,T}) \int_0^TK(s)\d s}\\
  &\le \big\{\psi_1(\|\mu\|_{1,T})\rr_\pp(x)\big\}^2,\ \ \ t\in [0,T].\end{align*} 
   Combining this with Jensen's inequality,     we prove 
 \eqref{*1N} with $\psi=\psi_1$ holds for $\rr_\pp(x)\ge \ff{r_0}2.$
   
(b)   Let $\rr_\pp(x)<\ff{r_0}2$.  
Simply denote $\aa=\aa(\|\mu\|_{1,T})$ and define 
$$\hat\tau:=\inf\{t\ge 0: \rr_\pp(X_t^\mu)\ge r_0\}. $$ By $(B_2)$ and It\^o's formula, 
we obtain
$$\d\rr_\pp(X_t^\mu)\le \aa \d t+ \d M_t,\ \ \ t\in [0,T\land \tau(X^\mu)\land\hat\tau]$$
for some martingale satisfying \eqref{*N}.  So,
$$ \E^x[\rr_\pp(X_{t\land\tau(X^\mu)\land\hat\tau}^\mu)] \le\aa\, \E^x[ t\land\tau(X^\mu)\land\hat\tau].$$
Combining this with step (a) and the strong Markov property, we obtain
  \beg{align*}  &\E^x[\rr_\pp(X_t^\mu)]= \E^x[\rr_\pp(X_{t\land\tau(X^\mu)}^\mu)]\le \E^x[\rr_\pp(X_{t\land\tau(X^\mu)\land\hat\tau}^\mu)]+ \E^x\Big[\1_{\{t\land\tau(X^\mu)\ge\hat\tau\}} \GG_{t-\hat\tau}(X_{\hat\tau}^\mu)\Big]\\
&\le \aa\, \E^x[ t\land\tau(X^\mu)\land\hat\tau]+ \P^x\big(t\land\tau(X^\mu)\ge\hat\tau\big) \psi_1(\|\mu\|_{1,T})r_0.
\end{align*} 
 Combining this with \eqref{LOT1} and \eqref{B*2}, we prove \eqref{*1N} for some increasing function $\psi: [0,\infty)\to (0,\infty),$  which is bounded if so is $\aa$. 
     \end{proof} 
 
 \beg{proof}[Proof of Theorem \ref{T2.2}] Let $X_t$ solve \eqref{E1} for $O$-distributions in $\scr P_O^1$. Then $X_t=X_t^\mu$ for $\mu_t:=\L_{X_t}^O$, so that 
 $$\|\mu_s\|_1= \E[\1_{O}(X_s)| X_s|]= \E[\1_{\{t<\tau(X)\}}|X_s|],\ \ s\in [0,T].$$ Combining this with \eqref{MES'}, we obtain  
\beg{align*} \|\mu_t\|_1^2&\le \Big(\E\ss{\E[\1_{\{t<\tau(X)\}}|X_t|^2|\F_0]}\Big)^2 \le 2 \Big(\E\ss{c (1+ \1_O(X_0)|X_0|^2)} \Big)^2+  2 c\int_0^t K(s) \|\mu_s\|_1^2\d s\\
&\le 2 c (1+ \E[\1_O(X_0)|X_0|])^2 + 2 c\int_0^t K(s) \|\mu_s\|_1^2\d s,\ \ t\in [0,T].\end{align*}
 By Gronwall's inequality, we find a constant $c_1>0$ such that  
 $$\sup_{t\in [0,T]}  \|\mu_t\|_1^2\le c_1 (1+ \E[\1_O(X_0)|X_0|])^2.$$ 
This together  \eqref{MES'} yields    \eqref{A01'} for some different constant $c>0$. It remains to prove the well-posedness and \eqref{A02'}. 
 
 (a) Well-posedness.  Let $\gg:=\L_{X_0}^O\in \scr P_O^1$. For any $N>0$, let
\beq\label{CNG} \C_N^\gg:=\bigg\{\mu\in C([0,T];\scr P_O^1):\ \mu_0=\gg, \  \sup_{t\in [0,T]} \e^{-Nt} \|\mu_t\|_1\le N\bigg\}.\end{equation} 
 We first observe that for some constant $N_0>0$, 
 \beq\label{IN} \Phi \C_N^\gg\subset \C_N^\gg,\ \ N\ge N_0.\end{equation}  Let $\mu\in \C_N^\gg$ and let $X_t^\mu$ solve \eqref{EM} for $X_0^\mu=X_0$. 
 Then $(\Phi\mu)_t=\L_{X_t^\mu}^O$. By \eqref{MES'} and
$$\|(\Phi\mu)_t\|_1 \le \E \ss{\E[\1_{O}(X_0)|X_{t\land\tau(X^\mu)}|^2|\F_0]},$$   we find a constant $c_1>0$ such that   
 $$\|(\Phi\mu)_t\|_1 \le      c_1  (1+\|\gg\|_1) +c_1 \bigg(  \int_0^t \|\mu_s\|_1^2\d s\bigg)^{\ff 1 2},\ \ t\in [0,T].$$
 Then for any $N\ge N_0:= c_1+ 2 c_1(1+\|\gg\|_1),$ we have 
 \beg{align*} &\sup_{t\in [0,T]}\e^{-Nt} \|(\Phi\mu)_t\|_1 \le c_1 (1+\|\gg\|_1)+ c_1\sup_{t\in [0,T]}\bigg(\int_0^t \e^{-2Ns}\|\mu_s\|_1^2 \e^{-2N(t-s)}\d s\bigg)^{\ff 1 2} \\
 &\le  c_1(1+\|\gg\|_1)+ c_1 N\sup_{t\in [0,T]}\bigg(\int_0^t \e^{-2N(t-s)}\d s\bigg)^{\ff 1 2} \\
 &\le c_1(1+\|\gg\|_1)+ c_1\ss N\le N.\end{align*}  
 Next, for any $N\ge N_0$, we intend to prove that $\Phi$ is contractive  in $\C_N^\gg$ under the  complete metric
 $$\W_{1,\theta}(\mu^1,\mu^2):=\sup_{t\in [0,T]}\e^{-\theta t} \W_1(\mu_t^1,\mu_t^2)$$
 for large enough $\theta>0$, so that  $\Phi$ has a unique fixed point in $\C^\gg=\cup_{N\ge N_0}\C_N^\gg$, hence  the well-posedness follows from Theorem \ref{T1}. 
 
 To this end, let $\mu^i\in \C_N^\gg$ and $X_t^i$ solve \eqref{EM} for $\mu=\mu^i$ and $X_0^i=X_0, i=1,2. $ By Lemma \ref{L2.2'} and noting that
 $\rr_\pp(x)\le |x-y|$ for $x\in O$ and $y\in\pp O$, we find a constant $c_2>0$ depending on $N$ such that for any $\mu^1,\mu^2\in \C_N^\gg$,
\beg{align*} &\E\big[\rr_\pp(X_t^1)\1_{\{t\land\tau_1\ge \tau_2\}}+ \rr_\pp(X_t^2)\1_{\{t\land \tau_2\ge \tau_1\}}\big]\\
 &\le c_2 \E\big[\rr_\pp(X_{t\land\tau_{1,2}}^1)\1_{\{t\land\tau_1\ge \tau_2\}}+ \rr_\pp(X_{t\land\tau_{1,2}}^2)\1_{\{t\land \tau_2\ge \tau_1\}}\big]\le 2 c_2 \E[|X_{t\land\tau_{1,2}}^1-X_{t\land\tau_{1,2}}^2|].\end{align*}
 Combining this with \eqref{C1'} and \eqref{PW'}, we find a constant $c_3>0$ depending on $N$ such that
\beq\label{WW1} \W_1((\Phi\mu^1)_t,(\Phi\mu^2)_t)\le c_3 \E[|X_0^1-X_0^2|] + c_3\bigg(\int_0^t K(s) \W_1(\mu_s^1,\mu_s^2)^2\d s\bigg)^{\ff 1 2},\ \ \mu^1,\mu^2\in \C_N^\gg.\end{equation} Since    $X_0^1=X_0^2=X_0$, 
 this implies the contraction of $\Phi$ in $\W_{1,\theta}$ for large enough $\theta>0.$

(b) Estimate \eqref{A02'}.  Now,  for $\mu_0^1,\mu_0^2\in \scr P_O^1$, let  $X_0^1,X_0^2$ be $\F_0$-measurable random variables on $\bar O$ such that
\beq\label{OLL}\L_{X_0^1}^O=\mu_0^1, \ \ \L_{X_0^2}^O=\mu_0^2,\ \ \E[|X_0^1-X_0^2|]= \W_1(\mu_0^1,\mu_0^2).\end{equation} 
Letting $X_t^i$ solve \eqref{E1} with initial value $X_0^i$, then $\mu^i:= (P_t^{D*}\mu_0^i)_{t\in [0,T]}$ is the unique  fixed point of $\Phi$ in $\C^{\mu_0^i},$ so that 
\beq\label{BGG} \mu_t^i=\L_{X_t^i}^O=\Phi \mu_t^i= P_t^{D*} \mu_0^i,\ \ i=1,2, t\in [0,T].\end{equation} 
When  $\aa$ is bounded,  \eqref{WW1} holds for some     constant $c_3>0$ independent of $N$, which together with \eqref{OLL} yields 
\beg{align*} &\W_1(\mu_t^1,\mu_t^2)= \W_1((\Phi\mu^1)_t,(\Phi\mu^2)_t) \le c_3 \E[|X_0^1-X_0^2|] + c_3\bigg(\int_0^t K(s) \W_1(\mu_s^1,\mu_s^2)^2\d s\bigg)^{\ff 1 2}\\
&= c_3 \W_1(\mu_0^1,\mu_0^2)+ c_3\bigg(\int_0^t K(s) \W_1(\mu_s^1,\mu_s^2)^2\d s\bigg)^{\ff 1 2},\ \ t\in [0,T].\end{align*}
By Gronwall's inequality and \eqref{BGG}, we obtain 
$$ \W_1(P_t^{D*}\mu_0^1, P_t^{D*}\mu_0^2)^2=\W_1(\mu^1_t,\mu^2_t)^2\le 2 c_3^2  \W_1(\mu_0^1,\mu_0^2)^2 \e^{2c_3^2 \int_0^t K(s)\d s},\ \ t\in [0,T].$$
Then the proof is finished. 
 \end{proof}

 \section{Singular case with distribution dependent noise} 
 
 In this part, we assume that $\si$ and $b$ are extended to $[0,T]\times\R^d\times \scr P_O$ but may be singular in the space variable. 
 To measure the   singularity, we recall locally integrable functional spaces introduced in \cite{XXZZ}.
For any $t>s\ge 0$ and $p,q\in (1,\infty)$, we write $f\in  \tilde{L}_p^q([s,t])$  if $f: [s,t]\times\R^d\to \R$ is measurable with
 $$\|f\|_{\tilde{L}_p^q([s,t])}:= \sup_{z\in\R^d}\bigg\{\int_{s}^t \bigg(\int_{B(z,1)}|f(u,x)|^p\d x\bigg)^{\ff q p} \d u\bigg\}^{\ff 1q}<\infty,$$ where $B(z,1):=\{x\in \R^d: |x-z|\le 1\}$ is the unit ball centered at point $z$.
 When $s=0$, we simply denote
\beq\label{KK0} \tt L_p^q(t)=\tt L_p^q([0,t]),\ \ \|f\|_{\tilde{L}_p^q(t)}=\|f\|_{\tilde{L}_p^q([0,t])}.\end{equation}
 We will take $(p,q)$ from the   space
 \beq\label{KK1} \scr K:=\Big\{(p,q): p,q>2, \ff{d}p+\ff 2 q<1\Big\}.\end{equation}
 For any $\mu\in C([0,T];\scr P_O)$,  let
\beq\label{SB}\si_t^\mu(x):= \si_t(x,\mu_t),\ \ \ b_t^\mu(x):= b_t(x,\mu_t)= b_t^{\mu,0}(x)+ b_t^{(1)}(x),\ \ \ (t,x)\in [0,T]\times \R^d,\end{equation}
 where $b^{\mu,0}_t(\cdot)$ is singular  and $b^{(1)}_t(\cdot)$ is Lipschitz continuous.

 As in the last section, we consider \eqref{E1} for $O$-distributions in $\scr P_O$ and $\scr P_O^1$ respectively.
 
 \subsection{For $O$-distributions in $\scr P_O$}

 \beg{enumerate} \item[{\bf (C)}] There exist    $K\in (0,\infty)$, $l\in\mathbb N$, $\{(p_i,q_i): 0\le i\le l\}\subset \scr K$ and
$1\le f_i\in \tt L^{q_i}_{p_i}(T)$ for $0\le i\le l,$ such that $\si^\mu$ and $b^\mu$ in \eqref{SB} satisfy the following conditions. 
 \item[$(C_1)$] For   any $\mu\in C([0,T];\scr P_O)$,   $a^\mu:= \si^\mu(\si^\mu)^*$ is invertible  with
$\|a^\mu\|_\infty +\|(a^\mu)^{-1}\|_\infty\le K$ and
$$\lim_{\vv\downarrow 0}\sup_{\mu\in C([0,T];\scr P_O)} \sup_{t\in [0,T], |x-y|\le \vv} \|a^\mu_t(x)-a^\mu_t(y)\|= 0.$$
\item[$(C_2)$] $b^{(1)}(0)$ is  bounded on $[0,T]$,   $\si^\mu_t$ is weakly differentiable for   $\mu\in C([0,T];\scr P_O)$,   and 
\beg{align*}& |b_t^{\mu,0}(x)|\le f_0(t,x),\ \
  \|\nn \si_t^\mu(x)\|\le \sum_{i=1}^l f_i(t,x),\  \\
  &|b_t^{(1)}(x)- b_t^{(1)}(y)|\le K |x-y|,\ \ t\in [0,T], x,y\in \R^d.\end{align*} 
  \item[$(C_3)$]  For any $t\in [0,T], x\in \R^d$ and $\mu,\nu\in \scr P_O$,
$$\|\si_t(x,\mu)-\si_t(x,\nu)\|+ |b_t(x,\mu)-b_t(x,\nu)|\le  \hat \W_1(\mu,\nu)\sum_{i=0}^l f_i(t,x).$$
  \end{enumerate}

 \beg{thm}\label{T3.1} Assume   {\bf (C)} and  $(A_2)$. Then the following assertions hold. 
 \beg{enumerate} \item[$(1)$] $\eqref{E1}$ is well-posed for $O$-distributions in $\scr P_O$.
\item[$(2)$]    For any $p\ge 1$, there exists a constant $c_p>0$ such that for any solution  $X_t$  to $\eqref{E1}$ for $O$-distributions in $\scr P_O$,  
 \beq\label{A1} \E \bigg[\sup_{t\in [0,T]} |X_t|^p \bigg|\F_0\bigg]= \E \bigg[\sup_{t\in [0,T]} |X_{t\land\tau(X)}|^p \bigg|\F_0\bigg]\le c_p\big(1+|X_0|^p\big).\end{equation} 
 \item[$(3)$] There exists a constant $c>0$ such that $\eqref{A02}$ holds.   \end{enumerate} \end{thm}
 
 For any $\mu\in C([0,T]; \scr P_O)$, instead of \eqref{EM} we consider the following SDE on $\R^d$:
 \beq\label{EM'} \d X_t^\mu= b_t^\mu(X_t^\mu)\d t+ \si_t^\mu(X_t^\mu)\d W_t,\ \ t\in [0,T].\end{equation} 
Noting that $\tt X_t^\mu:=X_{t\land\tau(X^\mu)}^\mu$ solves \eqref{EM},   the map $\Phi$ in \eqref{SB0}  is given by
$$ (\Phi \mu)_t:= \L_{X^\mu_{t\land\tau(X^\mu)}}^O,\ \ t\in [0,T].$$ So,     \eqref{PW}  and \eqref{PW'} remain true for $X_t^i$ solving  \eqref{EM'} with $\mu=\mu^i\in C([0,T];\scr P_O), i=1,2.$ 

By  \cite[Theorem 2.1]{21Ren}, see also \cite[Theorem 1.1]{W21e} for the distribution  dependent setting,   $(C_1)$ and $(C_2)$  imply that   this SDE is well-posed, and for any $p\ge 1$ there  exists a constant $c_p>0$ such that 
  \beq\label{MESN} \E \bigg[\sup_{t\in [0,T]} |X_t^\mu|^p \bigg|\F_0\bigg]\le c_p\big(1+|X_0^\mu|^p\big),\ \ \mu\in C([0,T];\scr P_O).\end{equation} 
  We have  the following lemma.
 
 \beg{lem}\label{L3.1} Assume {\bf (C)}. Then for any $j\ge 1$ there exists a constant $c>0$ and a function 
 $\vv: [1,\infty)\to (0,\infty)$ with $\vv(\theta)\downarrow 0$ as $\theta\uparrow \infty$, such that for any $\mu^1,\mu^2\in C([0,T]; \scr P_O)$ and any $X_t^i$ solving 
 $\eqref{EM'}$ with $\mu=\mu^i, i=1,2$, 
 $$\E\Big[\sup_{s\in [0,t]} |X_s^1-X_s^2|^j\big|\F_0\Big]\le c |X_0^1-X_0^2|^j+ \vv(\theta) \e^{j\theta t} \hat \W_{1,\theta} (\mu^1,\mu^2)^j,\ \ \theta\ge 1.$$  \end{lem}
 
 \beg{proof} The assertions follows from the proof of \cite[Lemma 2.1]{HW21c} for $\mu^i=\nu^i$ and for $\hat \W_1$ replacing $\W_{k}$ and $\W_{k,var}$.
 We figure it out for completeness. 
 
By \cite[Theorem 2.1]{YZ}, $(C_1)$ and $(C_2)$ imply that for large enough $\ll\ge 1$,  the  PDE 
  \beq\label{W*2}
\Big(\pp_t +\ff 12 {\rm tr}\{a_t^{\nu^1}\nn^2\}+b_t^{\mu^1} \cdot\nn  \Big) u_t  =  \ll u_t- b_t^{\mu^1,0}, \ \ t\in [0,T], u_T=0\end{equation} for $u: [0,T]\times\R^d\to \R^d$ 
has a unique solution such that
   \beq\label{3AA} \|\nn^2u\|_{\tt L_{p_0}^{q_0}(T)} \le c_0,\ \ \|u\|_\infty+\|\nn u\|_\infty\le \ff 1 2.\end{equation}
Let $Y_t^i:= \Theta_t(X_t^{i}), i=1,2, \Theta_t:=id +u_t$.
By It\^o's formula we obtain
{  \beg{align*}
  &\d Y_t^1=\big\{b_t^{(1)} + \lambda u_t \big\} (X_{t}^{1}) \d t+ (\{\nabla\Theta_t\}\sigma_t^{\nu^1} )(X_{t}^{1})\,\d W_t,\\
 & \d Y_t^2 =\big\{\big\{b_t^{(1)} + \lambda u_t+{  (\nabla\Theta_t)(b^{\mu^2}_t-b^{\mu^1}_t)\big\}(X_{t}^{2})} \\
 &  + {  \ff 1 2\big[{\rm tr} \{(a_t^{\nu^2}- a_t^{\nu^1}) \nn^2   u_t\}\big](X_{t}^{2})} \big\}\d t+(\{\nabla\Theta_t\}\sigma_t^{\nu^2})(X_{t}^{2})\,\d W_t.\end{align*}}
Let $\eta_{t}:= |X_{t}^{1}- X_{t}^{2}|$ and
   \beg{align*} &g_r:=  \sum_{i=0}^l f_i  (r, X_{r}^{2}),\ \ \tt g_r:= g_r  \|\nn^2  u_r(X_{r}^{2})\|,\\
   &\bar g_r:= \sum_{i=1}^2 \|\nn^2 u_r\|(X_r^i)+\sum_{j=1}^2 \sum_{i=0}^l f_i  (r, X_{r}^{j}),\ \ r\in [0,T].\end{align*}
 Since $b_t^{(1)} + \lambda u_t $ is Lipschitz continuous uniformly in $t\in [0,T]$, by  {\bf (C)}  and the maximal functional inequality in \cite[Lemma 2.1]{XXZZ},
     there exists a constant
 $c_1>0$  such that
   \beg{align*}&\big| \big\{b_r^{(1)} + \lambda u_r \big\} (X_{r}^{1})  - \big\{b_r^{(1)} + \lambda u_r \big\} (X_{r}^{2}) \big|\le c_1 \eta_r,\\
 & \big| \big\{(\nabla\Theta_r)(b^{\mu^2}_r-b^{\mu^1}_r)\big\}(X_{r}^{2})\big|\le c_1  g_r \hat \W_{1} (\mu^1_r,\mu^2_r),\\
 &\big|\big[{\rm tr} \{(a_r^{\nu^2}- a_r^{\nu^1}) \nn^2   u_r\}\big](X_{r}^{2})\big|\le c_1 \tt g_r \hat \W_{1}(\mu_r^1,\mu_r^2),\\
 &\big\| \big\{(\nabla\Theta_r ) \sigma_r^{\nu^1} \big\} (X_{r}^{1})- \big\{(\nabla\Theta_r ) \sigma_r^{\mu^2} \big\} (X_{r}^{2})\big\|\\
 &\le c_1  \bar g_r \eta_r+  c_1 g_r \hat \W_{1}(\mu^1_r,\mu^2_r),\ \ r\in [0,T].\end{align*}
So, by It\^o's formula, for any $j\ge k$
   we find a constant    $c_2>1$     such that
\beq\label{**1}  \d |Y_t^1-Y_t^2|^{2j}\le c_2  \eta_t^{2j} \d A_t + c_2 (g_t^2+\tt g_t) \hat \W_{1}(\mu^1_t,\mu_t^2)^{2j}   \d t +\d M_t\end{equation}
   holds for some martingale $M_t$   with $M_0=0$ and
   $$A_t:=\int_0^t \big\{1+g_s^2+\tt g_s+ \bar g_s^2\big\}\d s,\ \ t\in [0,T].$$
   Since $\|\nn u\|_\infty\le \ff 1 2$ implies $|Y_t^1-Y_t^2|\ge \ff 1 2 \eta_t,$ this implies
     \beq\label{NNP} \beg{split}& \eta_{t}^{2j}  \le 2^{2j}M_t+ 2^{2j} \eta_0^{2j} +  2^{2j}c_2\int_0^t \eta_{r}^{2j} \d A_r \\
     &+ 2^{2j}c_2 \int_0^t  (g_s^2+\tt g_s)   \hat \W_{1}(\mu^1_s,\mu_s^2)^{2j}  \d s,\ \ t\in [0,T]  \end{split}\end{equation}
  for some constant $c_2>0.$  
 By   \eqref{3AA},  $f_i\in \tt L^{q_i}_{p_i}(T)$ for $(p_i,q_i)\in \scr K$,   Krylov's estimate  (see \cite[Theorem 3.1]{YZ}) which implies   
  Khasminskii's estimate (see \cite[Lemma 4.1(ii)]{XXZZ}), we find an increasing function
 $\psi: (0,\infty)\to (0,\infty)$ and a decreasing function $\vv: (0,\infty)\to (0,\infty)$ with $\vv(\theta) \downarrow 0$ as $\theta\uparrow \infty$, such that
   $$\E[\e^{r A_T}|\F_0]\le \psi(r),\ \ \ r>0,$$ 
   $$  \sup_{t\in [0,T]} \E\bigg( \int_0^t \e^{-2k \theta(t-r)}(g_r^2+\tt g_r)\d r \bigg|\F_0\bigg)\le \vv(\theta), \ \ \theta>0.$$ 
  By  the stochastic Gronwall inequality
and  the maximal inequality (see \cite{XXZZ}), we find   a   constant     $c_3>0$  depending on   $N$  such that    \eqref{NNP}   yields
 \beg{align*}   &\Big\{\E\Big(\sup_{s\in [0,t]}  \eta_{s}^j\Big|\F_0\Big)\Big\}^{2} \\
& \le c_3 \E \bigg(\eta_0^{2j}+ \int_0^t  (g_s^2+\tt g_s)  \hat \W_{1}(\mu^1_s,\mu_s^2)^{2j}  \d s\bigg|\F_0\bigg) \\
 &\le c_3  \eta_0^{2j}+c_3  \e^{2j\theta t}   \vv(\theta)      \hat \W_{1}(\mu^1,\mu^2)^{2j},\ \ t\in [0,T], \theta>0.\end{align*}
 This finishes the proof. 
  \end{proof}   
 
 \beg{proof}[Proof of Theorem \ref{T3.1}]  Let $X_t$ solve \eqref{E1}. We have $X_t=X_{t\land\tau(X^\mu)}^\mu$ for $X_t^\mu$ solving \eqref{EM'} with
 $$X_0^\mu=X_0,\ \ \mu_t:= \L_{X_t}^O,\ \ t\in [0,T]. $$
 So, \eqref{A1} follows from \eqref{MESN}. It remains  to prove the well-posedness and estimate \eqref{A02}. 
 
 (a) Well-posedness. Let  $X_0$ be an $\F_0$-measurable random variable on $\bar O$, and let $\C^\gg$ be in \eqref{CG} for $\gg=\L_{X_0}^O$. By Theorem \ref{T1}, it suffices to prove that $\Phi$ is contractive in $\C^\gg$ under $\hat \W_{1,\theta}$ for large enough $\theta>0.$ 
 
 By \eqref{PW}, \eqref{ORR} and Lemma \ref{L3.1} for $X_0^1=X_0^2=X_0,$ we find a constant $c_1>0$ such that
 $$\hat W_1((\Phi\mu^1)_t,(\Phi\mu^2)_t) \le c_1 \vv(\theta) \hat\W_1(\mu^1,\mu^2),\ \ \ \mu^1,\mu^2\in \C^\gg.$$
 Since $\vv(\theta)\to 0$ as $\theta\to\infty$, $\Phi$ is $\hat\W_{1,\theta}$-contractive for large enough $\theta>0.$
 
 (b) Estimate \eqref{A02}. Let $X_t^1,X_t^2$ solve \eqref{E1} with $X_0^1,X_0^2$ satisfying \eqref{X0}. Then $$(\Phi\mu^i)_t= \mu^i_t:=\L_{X^i}^O=P_t^{D*}\mu^i,\ \ \  i=1,2,$$
 so that  \eqref{PW}, \eqref{ORR} and Lemma \ref{L3.1} imply
 $$\hat W_1(\mu^1,\mu^2)= \hat W_1((\Phi\mu^1)_t,(\Phi\mu^2)_t) \le c_1 \hat \W_1(\mu_0^1,\mu_0^2)+ c_1 \vv(\theta) \hat\W_1(\mu^1,\mu^2),\ \ t\in [0,T]$$ for some constant $c_1>0$. 
 Taking $\theta>0$ large enough such that $\vv(\theta)\le \ff 1 {2 c_1}$, we derive \eqref{A02} for some constant $c>0.$

 \end{proof} 
  
   \subsection{For $O$-distributions in $\scr P_O^1$} 
   
    \beg{enumerate} \item[{\bf (D)}] There exist  an increasing function   $\aa: [0,\infty)\to (0,\infty)$, constants $K>0, l\in\mathbb N$, $\{(p_i,q_i): 0\le i\le l\}\subset \scr K$ and
functions $1\le f_i\in \tt L^{q_i}_{p_i}(T)$ for $0\le i\le l$ such that  $\si^\mu$ and $b^\mu$ in \eqref{SB} satisfy the following conditions. 
 \item[$(D_1)$] For   any $\mu\in C([0,T];\scr P_O^1)$,   $a^\mu:= \si^\mu(\si^\mu)^*$ is invertible  with
\beg{align*} &\|a^\mu\|_\infty +\|(a^\mu)^{-1}\|_\infty\le \aa(\|\mu\|_{1,T}),\\
&\lim_{\vv\downarrow 0}\sup_{\mu\in C([0,T];\scr P_O^1)} \sup_{t\in [0,T], |x-y|\le \vv} \|a^\mu_t(x)-a^\mu_t(y)\|= 0.\end{align*}
\item[$(D_2)$] $b^{(1)}(0)$ is  bounded on $[0,T]$,   $\si^\mu_t$ is weakly differentiable for  $\mu\in C([0,T];\scr P_O^1)$,   and 
\beg{align*}& |b_t^{\mu,0}(x)|\le f_0(t,x)+\aa(\|\mu\|_{1,T}),\ \
  \|\nn \si_t^\mu(x)\|\le \sum_{i=1}^l f_i(t,x)+\aa(\|\mu\|_{1,T}), \\
  &|b_t^{(1)}(x)- b_t^{(1)}(y)|\le K |x-y|,\ \ t\in [0,T], x,y\in \R^d.\end{align*} 
  \item[$(D_3)$]  For any $t\in [0,T], x\in \R^d$ and $\mu,\nu\in \scr P_O$,
$$\|\si_t(x,\mu)-\si_t(x,\nu)\|+ |b_t(x,\mu)-b_t(x,\nu)|\le  \W_1(\mu,\nu)\sum_{i=0}^l f_i(t,x).$$
 \item[$(D_4)$] There exists $r_0\in (0,1]$ such that  $\rr_\pp \in C_b^2(\pp_{r_0}O)$, and  for any $\mu\in C([0,T];\scr P_O^1)$, 
\beq\label{LU'}    \<b_t^\mu(x),\nn  \rr_\pp(x)\>\le   \aa(\|\mu\|_1),\ \ \ x\in \pp_{r_0}O,  \end{equation} 
\beq\label{LU2'}
 \<b_t^\mu(x), x-y\>\le \aa(\|\mu\|_{1,T})(f_0(t,x)^2+|x-y|^2),\ \ x\in O, y\in\pp O, t\in [0,T].\end{equation}
  \end{enumerate}
Note that when $b^{(1)}=0$, \eqref{LU'}  is implied by the first condition in $(D_2)$. 
   
   \beg{thm}\label{T3.2} Assume   {\bf (D)}. Then the following assertions hold. 
 \beg{enumerate} \item[$(1)$] $\eqref{E1}$ is well-posed for $O$-distributions in $\scr P_O^1$.
\item[$(2)$]    For any $p\ge 1$, there exists a constant $c_p>0$ such that for any solution  $X_t$  to $\eqref{E1}$ for $O$-distributions in $\scr P_O$,  
 \beq\label{A31} \E \bigg[\sup_{t\in [0,T]} |X_t|^p \bigg|\F_0\bigg]\le c_p\big\{1+|X_0|^p+ \big(\E[\1_O(X_0)|X_0|]\big)^p\big\}.\end{equation} 
 \item[$(3)$] If $\aa$ is bounded, then there exists a constant $c>0$ such that $\eqref{A02'}$ holds.   \end{enumerate} \end{thm}

 By  the proof of \cite[(2.17)]{HW21c},    {\bf (D)} implies that for any  $\mu\in C([0,T]; \scr P_O^1)$, 
 the SDE \eqref{EM'}    is well-posed, and for any $p\ge 1$ there 
    exists a constant $c_p>0$ such that 
  \beq\label{MESN'} \E \bigg[\sup_{t\in [0,T]} |X_t^\mu|^{2p} \bigg|\F_0\bigg]\le c_p\bigg\{1+|X_0^\mu|^{2p}+\int_0^T\|\mu_s\|_1^{2p}\d s\bigg\},\ \ \mu\in C([0,T];\scr P_O).\end{equation} 
 For any $\mu^1,\mu^2\in \scr P_O^1$, let $X_t^i$ solve \eqref{EM'} for $\mu=\mu^i, i=1,2.$ 
  
For any $N>0$ and $\gg\in \scr P_O^1$, let $\C_N^\gg$ be in \eqref{CNG}.    Since restricting to $\mu,\nu\in\C_N^\gg$ the conditions in {\bf (D)} hold for 
a constant $\aa_N$ replacing the function $\aa$, by  repeating the proof of Lemma \ref{L3.1} with $\W$ replacing $\hat \W$,  we prove that the following  result.

\beg{lem}\label{L3.2} Assume {\bf (D)}. 
For any $N>0$ and  $j\ge 1$, there exists a constant $c>0$ and a function 
 $\vv: [1,\infty)\to (0,\infty)$ with $\vv(\theta)\downarrow 0$ as $\theta\uparrow \infty$, such that for any $\mu^1,\mu^2\in \C_N^\gg$ and any $X_t^i$ solving 
 $\eqref{EM'}$ with $\mu=\mu^i, i=1,2$, 
 $$\E\Big[\sup_{s\in [0,t]} |X_s^1-X_s^2|^j\big|\F_0\Big]\le c |X_0^1-X_0^2|^j+ \vv(\theta) \e^{j\theta t} \hat \W_{1,\theta} (\mu^1,\mu^2)^j,\ \ \theta\ge 1.$$ 
 When $\aa$ is bounded, the constant $c$ does not depend on $N$. \end{lem} 
  Moreover, we need  the following  result analogous to Lemma \ref{L2.2'}.
  
  \beg{lem}\label{L3.3} Assume {\bf (D)}. Then the assertion in Lemma \ref{L2.2'} holds. \end{lem} 
 
\beg{proof} It suffices to prove \eqref{*1N} for some increasing function $\psi$ which is bounded if so is $\aa$.

(a) Let $\rr_\pp(x)\ge \ff{r_0}2$ and $y\in\pp O$ such that $\rr_\pp(x)=|y-x|. $
 By \eqref{LU2'} and $(D_2)$,    we find an increasing function $\psi_1: [0,\infty)\to (0,\infty)$ which is bounded if so is $\aa$, such that  
 $$\d |X_t^\mu-y|^2\le  \psi_1(\|\mu\|_{1,T}) \Big(\sum_{i=0}^l f_i(t,X_t^\mu)^2+|X_t^\mu-y|^2 \Big)\d t +\d M_t,\ \ t\in [0,T\land\tau(X^\mu)]$$
   for some martingale $M_t$.  Next, by  \cite[Theorem 3.1]{YZ}, {\bf (D)} implies that for some  increasing function $\psi_2: [0,\infty)\to (0,\infty)$ which is bounded if so is $\aa$,  
   the following Krylov's estimate holds:
 $$\E\bigg(\int_0^T f_i(t,X_t^\mu)^2\d t\bigg|\F_0\bigg)\le \psi_2(\|\mu\|_{1,T}) \|f_i\|_{\tt L_{q_i}^{p_i}(T)}^2,\ \ 0\le i\le l.$$
   Combining these with  $|x-y|=\rr_\pp(x)$,  we derive 
\beg{align*}  \E[|X_t^\mu-y|^2|\F_0]\le&\,   \rr_\pp(x)^2+ \psi_1(\|\mu\|_{1,T}) \psi_2(\|\mu\|_{1,T}) \sum_{i=0}^l \|f_i\|_{\tt L_{q_i}^{p_i}(T)}^2 \\
&+  \psi_1(\|\mu\|_{1,T})  \int_0^t \E[|X_s^\mu-y|^2|\F_0]\d s,\ \ t\in [0,T].\end{align*} 
 By Gronwall's inequality and $\rr_\pp(x)\ge \ff{r_0}2$,   we find an  increasing function $\psi: [0,\infty)\to (0,\infty)$ which is bounded if so is $\aa$,  such that
 $$\E[|X_t^\mu-y|^2|\F_0]\le \psi(\|\mu\|_{1,T}) \rr_\pp(x).$$
Since   $\rr_\pp(X_t^\mu)\le |X_t^\mu-y|$,  we prove   \eqref{*1N}  for $\rr_\pp(x)\ge \ff{r_0}2. $
   
(b)   Let $\rr_\pp(x)<\ff{r_0}2$.  By   $(D_1)$,  \eqref{LU'} and $\rr_\pp\in C_b^2(\pp_{r_0}O)$, \eqref{LU} holds for some different increasing function $\aa$ which is bounded if so is the original one. Then  step $(b)$ in  proof of Lemma \ref{L2.2'}  implies  the desired estimate. 
 \end{proof} 

\beg{proof}[Proof of Theorem \ref{T3.2}]  Let $X_t$ solve \eqref{E1} for $O$-distributions in $\scr P_O^1$. We have $X_t=X_{t\land\tau(X^\mu)}^\mu$ for $X_t^\mu$ solving 
\eqref{EM'} with
$$X_0^\mu=X_0,\ \ \mu_t= \L_{X_t}^O,\ \ t\in [0,T].$$
So,  as explained in the beginning of proof of Theorem \ref{T2.2} that \eqref{A31} follows from \eqref{MESN'}. It suffices to prove the well-posedness and estimate \eqref{A02'}. 

   Let $X_0$ be an $\F_0$-measurable random variable with $\L_{X_0}^O\in \scr P_O^1$, and let $\C_N^\gg$ be in \eqref{CNG} for $N>0.$   
By the proof of \cite[Lemma 2.2(1)]{HW21c}, there exists $N_0>0$ such that  $\Phi \C_N^\gg\subset \C_N^\gg$ for any $N\ge N_0$. For the well-posedness,  it suffices to prove that for any $N\ge N_0$, 
$\Phi$ is contractive in $\C_N^\gg$ under the metric $\W_{1,\theta}$ for large enough $\theta>0$.  This follows from  \eqref{PW'}, Lemma \ref{L3.2} and Lemma \ref{L3.3}.

Finally, by  using $\W_1$ replacing $\hat \W_1$ in step (b) in the proof of Theorem \ref{T3.1},   \eqref{A02'} follows from Lemma \ref{L3.2} with $c$ independent of $N$. 

\end{proof}

\section{Singular case with distribution independent noise }

In this part, we let $\si_t(x,\mu)=\si_t(x)$ do not depend on $\mu$, so that \eqref{E1} becomes
\beq\label{E2} \d X_t =\1_{\{t<\tau(X)\}}\big\{b_t(X_t, \L_{X_t}^O)\d t +\si_t(X_t)\d W_t\big\},\ \ t\in [0,T].\end{equation}
In this case,  we are able to study the well-posedness of the equation on  an arbitrary connected open  domain $O$, for which  we only need  
$b_t(x,\cdot)$ to be   Lipschitz continuous  with respect to a weighted variation distance.
  
 For  a measurable function $V: O\to  [1,\infty)$, let 
 $$\scr P_O^V:= \bigg\{\mu\in \scr P_O:\ \mu(V):=\int_O V\d\mu<\infty\bigg\}.$$
 This is a Polish space under the weighted variation distance 
\beq\label{WV} \|\mu-\nu\|_{V}:= \sup_{|f|\le V} |\mu(f)-\nu(f)|,\ \ \mu,\nu\in \scr P_O^V.\end{equation}
When $V\equiv 1$, $\|\cdot\|_V$  reduces to the total variation norm. 
 We  will take $V$ from the   class $\scr V$ defined as follows.

\beg{defn} We denote $V\in\scr V$, if $1\le V\in C^2(\R^d)$ such that the level set  $\{V\le r\}$ for $ r>0$ is compact, and there exist constants $K,\vv>0$ such that for any $x \in   O,$
 $$\sup_{y\in B(x,\vv)}\big\{|\nn V(y)|+   \|\nn^2 V(y)\| \big\}\le K V(x),$$   where $B(x,\vv):=\{y\in\R^d: |y-x|<\vv\}.$ \end{defn}

  \subsection{Main result} 

 \beg{enumerate}
\item[{\bf (E)}]     $\si$ has an extension to $[0,T]\times \R^d$ which is weakly differentiable in $x\in\R^d$, and $b$ has a decomposition $b_t(x,\mu)=b^{(0)}_t(x)+b^{(1)}_t(x,\mu)$,  such that  the following conditions hold.
\item[$(E_1)$]  $a:= \si\si^*$ is invertible with $\|a\|_\infty+\|a^{-1}\|_\infty<\infty$    and
$$ \lim_{\vv\to 0} \sup_{|x-y|\le \vv, t\in [0,T]} \|a_t(x)-a_t(y)\|=0.$$

 \item[$(E_2)$]   there exist $l\in \mathbb N$ and $1\le f_i \in \tt L_{q_i}^{p_i}(T)$ with $  (p_i,q_i) \in  \scr K, 0\le i\le l,$ such that 
$$ |\1_O b^{(0)}|\le f_0,\ \ \ \|\nn \si\|\le \sum_{i=1}^l f_i.$$
  \item[$(E_3)$]   there exists  $V\in \scr V$ such that  for any $\mu\in \C_V:=  C([0,T];\scr P_O^V)$,  $\1_O(x) b^{(1)}_t(x,\mu_t)$   is  locally bounded in $(t,x)\in [0,T]\times\R^d$. Moreover,  there exist constants $K,\vv>0$    such that  
\beg{align*}  
&   \<b^{(1)}(x,\mu),\nn V(x) \>+\vv |b^{(1)}(x,\mu)|\sup_{B(x,\vv)} \big\{|\nn V|+|\nn^2 V|\big\}\\
 &\le K \big\{V(x)+ \mu(V)\big\}  \ \ x\in O,\mu\in \scr P_O^V. \end{align*}
\item[$(E_4)$]   there exists a constant $\kk>0$ such that
 \beq\label{KK1}  \sup_{x\in O} |b_t(x,\mu)- b_t(x,\nu)|\le \kk \|\mu-\nu\|_V,\ \ \mu,\nu\in \scr P_O^V.\end{equation}
  \end{enumerate}

  \beg{thm}\label{T5.1}  Assume {\bf (E)}.
 Then   $\eqref{E2}$ is well-posed for $O$-distributions in
  $\scr P_O^V,$ and  for any $p\ge 1$, there exists a constant $c_p>0$ such that  any solution $X_t$ of $\eqref{E2}$ for  $O$-distributions in  $\scr P_O^V$ satisfies  
 \beq\label{KK2} \E\Big[\sup_{t\in [0,T]} V(X_t)^p\Big|\F_0\Big]\le c_p\,  V(X_0)^p. \end{equation}
   \end{thm}
\beg{proof} Let $\scr P^V(\bar O)$ be the space of all probability measures $\mu$ on $\bar O$ with $\mu(V)<\infty$, which is a Polish space under  the weighted variation distance defined in \eqref{WV} for $\mu,\nu\in \scr P^V(\bar O)$. 
  We extend 
   $b_t(x,\cdot)$ from $\scr P_O^V$   to  $\scr P^{V}(\bar O)$ by setting
 $$b_t(x,\mu):=b_t(x,  \mu(O\cap \cdot)),\ \ \mu\in \scr P^V(\bar O).$$
   Then {\bf (E)} implies the same assumption for $\scr P^V(\bar O)$ replacing $\scr P_O^V$. So, the desired assertions follow from Theorem \ref{T5.2} presented in the next subsection.
\end{proof} 

\subsection{An extension of Theorem \ref{T5.1}}
Consider the following SDE on $\bar O$: 
\beq\label{E2'} \d X_t =\1_{\{t<\tau(X)\}}\big\{b_t(X_t, \L_{X_t})\d t +\si_t(X_t)\d W_t\big\},\ \ t\in [0,T],\end{equation}
where $\tau(X):= \inf\{t\ge 0: X_t\in \pp O\}$ as before, and $\L_{X_t}$ is the distribution of $X_t$. 

The strong/weak solution of \eqref{E2'}  is defined as in Definition \ref{D1.1} with $\L$ replacing $\L^O$. 
We call this equation well-posed for distributions in $ \scr P^V(\bar O)$, if for any $\F_0$-measurable random variable $X_0$ on $\bar O$ with $\L_{X_0}\in \scr P^V(\bar O)$
(respectively, any $\mu_0\in \scr P^V(\bar O)$), \eqref{E2'} has a unique solution starting at $X_0$ (respectively, a unique weak solution with initial distribution $\mu_0$) such that
$\L_{X}=(\L_{X_t})_{t\in [0,T]}\in C([0,T]; \scr P^V(\bar O))$. 


  \beg{thm}\label{T5.2} Assume that {\bf (E)} holds for $\scr P^{V}(\bar O)$ replacing $\scr P_O^V$. Then
 $\eqref{E2'}$ is well-posed for distributions in $\scr P^V(\bar O)$ and \eqref{KK2} holds. \end{thm}

 \beg{proof}

  (1) Let $X_0$ be an    $\F_0$-measurable random variable on $\bar O$  with
  $$\gg:=\L_{X_0}\in \scr P_{V}(\bar O).$$
  Let
  $$\scr C_V^\gg(\bar O):=\big\{\mu\in C([0,T]; \scr P_{V}(\bar O)):\ \mu_0=\gg  \big\}.$$
For any $\mu\in \scr C_V^\gg(\bar O)$, let $X_t^\mu$ solve \eqref{E2} with $X_0^\mu=X_0$, i.e.  
  \beq\label{P1} \d X_t^\mu= \1_{\{t<\tau(X^\mu)\}} \big\{b_t(X_t^\mu,\mu_t)\d t+\si_t(X_t^\mu)\d W_t\big\},\ \ X_0^\mu=X_0, t\in [0,T].\end{equation}
Let $(\Phi\mu)_t:=\L_{X_t^\mu}, t\in [0,T].$  Then it suffices to prove that $\Phi$ has a unique fixed point in $\scr C_V^\gg(\bar O)$.  To this end, for any $N\ge 1$, let
$$\scr C_{V,N}^\gg(\bar O):= \Big\{\mu\in \scr C_V^\gg(\bar O): \sup_{t\in [0,T]}\e^{-Nt}\mu_t(V)\le N\gg(V)\Big\}.$$
It suffices to find a constant $N_0>0$ such that for any $N\ge N_0$, $\Phi$ has a unique fixed point in $\scr C_{V,N}^\gg(\bar O).$ We finish the proof by two steps.

(a) The $\Phi$-invariance of $\scr C_{V,N}^\gg(\bar O)$ for large $N$. 
For any $\ll\ge 0$and $N\ge 1$, $\scr C_{V,N}^\gg(\bar O)$ is a complete space under the metric
  $$\rr_\ll(\mu,\nu):=\sup_{t\in [0,T]} \e^{-\ll t} \|\mu_t-\nu_t\|_V,\ \ \mu,\nu\in \scr C_{V,N}^\gg(\bar O).$$ 
  Let $\mu\in \scr C_{V,N}^\gg(\bar O).$ By \eqref{P1}, {\bf (E)} with $V\in \scr V$ and It\^o's formula, for any $p\ge 1$ we find a constant $c_1(p)>0$ such that
  $$\d V(X_t^\mu)^p \le \1_{\{t<\tau(X^\mu)\}} \big\{\d M_t + c_1\big\{V(X_t^\mu)^p+ \mu_t(V)^p\big\}\d t\big\},\ \ t\in [0,T],$$
  where $M_t$ is a martingale with
  $$\d\<M\>_t\le c_1 V(X_t^\mu)^p\d t.$$
  By using BDG's and Gronwall's inequality, we find a constant $c_2(p)>0$ such that
  \beq\label{G1}\beg{split}& \E\Big[\sup_{s\in [0,t]}V(X_s^\mu)^p\Big]=\E\Big[\sup_{s\in [0,t\land\tau(X^\mu)]}V(X_s^\mu)^p\Big]\\
  &\le c_2(p) V(X_0)^p+c_2(p)\int_0^t \mu_s(V)^p\d s,\ \ t\in [0,T].\end{split}\end{equation} 
  Consequently, for $p=1$ and $c_2=c_2(1)$ we derive 
  $$(\Phi\mu)_t(V)= \E[V(X_t^\mu)] \le c_2  \gg(V) + c_2 \bigg(\int_0^t \mu_s(V)^2\d s\bigg)^{\ff 1 2},$$
  so that by $\mu\in \scr C_{V,N}^\gg(\bar O)$ we obtain
  \beg{align*}&\sup_{t\in [0,T]}\e^{-Nt} (\Phi\mu)_t(V)\le c_2 \gg(V)+ c_2 \sup_{t\in [0,T]} \bigg(\int_0^t\e^{-2N(t-s)} N^2\gg(V)^2\d s\bigg)^{\ff 1 2}\\
  &c_2\big(1+ \ss N\big)\gg(V)\le N\gg(V)\end{align*}
  provided $N\ge N_0$ for a large enough constant $N_0\ge 1.$  By the continuity of $X_t^\mu$ in $t$, $(\Phi\mu)_t$ is weakly continuous in $t$.  Therefore,  
  $$\Phi \scr C_{V,N}^\gg(\bar O)\subset \scr C_{V,N}^\gg(\bar O),\ \ \ N\ge N_0.$$ 
  
  (b) Let $N\ge N_0$. It remains to show that  $\Phi$ has a unique fixed point in $\scr C_{V,N}^\gg(\bar O).$  
  By \eqref{G1} with $p=2$ and $V\ge 1$, there exists a constant $c_3>0$  such that
  \beq\label{P2} \E\Big[\sup_{t\in [0,T]}V(X_t^\mu)^2\Big|\F_0\Big]\le c_3^2 V(X_0)^2,\ \ \mu\in \scr C_{V,N}^\gg(\bar O).\end{equation}
   For any $\mu^i\in \scr C_V(\bar O), i=1,2,$ we estimate $\|(\Phi\mu^1)_t-(\Phi\mu^2)_t\|_V$ by using
  Girsanov's  theorem. Let $X_t^1$ be the unique solution for the SDE
  \beq\label{R1}\d X_t^1 = \1_{\{t<\tau(X^1)\}}\big\{b_t(X_t^1, \mu_t^1)\d t +\si_t(X_t^1)\d W_t\big\},\ \ X_0^1=  X_0.\end{equation}
 By the definition of $\Phi$, we have
\beq\label{*P'} (\Phi\mu^1)_t= \L_{X_t^1},\ \ t\in [0,T].\end{equation} 
  To construct $(\Phi\mu^2)_t$ using Girsanov's theorem, let 
  $$\xi_t:= \1_{\{t<\tau(X^1)\}}\{\si_t^*(\si_t\si_t^*)^{-1}\}(X_t^1) \{b_t(X_t^1, \mu_t^2)-b_t(X_t^1, \mu_t^1)\},\ \ t\in [0,T].$$
  By {\bf (E)}, there exists a constant $k>0$ such that
\beq\label{xi} |\xi_t|\le k \|\mu_t^1-\mu_t^2\|_V,\ \ t\in [0,T].\end{equation}
  So, by Girsanov's theorem,
  $$\tt W_t:= W_t-\int_0^t \xi_s\d s,\ \ t\in [0,T]$$
  is an $m$-dimensional Brownian motion under the probability measure $\Q:= R_T\P$, where
  $$R_s:= \e^{\int_0^s \<\xi_t,\d W_t\>-\ff 1 2 \int_0^s|\xi_t|^2\d t},\ \ s\in [0,T].$$
  Reformulate \eqref{R1} as
  $$\d X_t^1 = \1_{\{t<\tau(X^1)\}}\big\{b_t(X_t^1, \mu_t^2)\d t +\si_t(X_t^1)\d \tt W_t\big\},\ \ X_0^1=\tt X_0.$$
  By the weak uniqueness of \eqref{P1}, we obtain 
  $$(\Phi\mu^2)_t= \Q(X_{t\land \tau(X^1)}^1\in \d x)= \L_{X_t^1|\Q}.$$
  Combining this with   \eqref{P2} and \eqref{*P'}, we derive
 \beq\label{PC}  \beg{split} & \|(\Phi\mu^1)_t-(\Phi\mu^2)_t\|_V\le \E\big[V(X_t^1)|R_t-1|\big]\\
&  \le \E\big[\{\E(V(X_t^1)^2|\F_0)\}^{\ff 1 2} \{\E(|R_t-1|^2|\F_0)\}^{\ff 1 2}\big]\\
  &\le c_3  \E\big[V(X_0)\{\E(|R_t-1|^2|\F_0)\}^{\ff 1 2}\big].\end{split}\end{equation} 
  On the other hand, by $\mu^1,\mu^2\in \scr C_{V,N}^\gg(\bar O)$, \eqref{xi}, and noting that $\e^r-1\le r\e^r$ for $r\ge 0$, we find a constant
  $c>0$ such that
  \beg{align*}&\E[|R_t-1|^2|\F_0]= \E[\e^{2\int_0^t \<\xi_s,\d W_s\>-  \int_0^t|\xi_s|^2\d s}-1|\F_0]\\
  &\le \E\big[\e^{2\int_0^t \<\xi_s,\d W_s\>- 2\int_0^t|\xi_s|^2\d t}|\F_0]\e^{k^2\int_0^t\|\mu_s^1
  -\mu_s^2\|_V^2\d s}-1\\
  &= \e^{k^2\int_0^t\|\mu_s^1 -\mu_s^2\|_V^2\d s}-1\le \e^{k^2\int_0^t\|\mu_s^1 -\mu_s^2\|_V^2\d s}\int_0^t k^2\|\mu_s^1-\mu_s^2\|_V^2\d s\\
  &\le c^2 \int_0^t  \|\mu_s^1-\mu_s^2\|_V^2\d s,\ \ \ t\in [0,T].\end{align*}
  Combining this with \eqref{PC} and  letting $C=cc_3\,\E[V(X_0)]$,  we arrive at
   \beg{align*} &\rr_\ll(\Phi(\mu^1),\Phi(\mu^2))\le C\sup_{t\in [0,T]}\e^{-\ll t} \bigg(\int_0^t\|\mu_s^1-\mu_s^2\|_V^2\d s\bigg)^{\ff 12}\\
   &\le C \rr_\ll(\mu^1,\mu^2) \bigg(\int_0^t \e^{-2\ll(t-s)}\d s\bigg)^{\ff 1 2}.\end{align*}
   Thus, when $\ll>0$ is large enough, $\Phi$ is contractive in $\rr_\ll$ and hence has a unique fixed point in $C_{V,N}^\gg(\bar O).$

   (3)  Uniqueness and  \eqref{KK2}.  It is easy to see that for any (weak) solution $X_t$ of \eqref{E2'} for distributions
   in $\scr P^V(\bar O)$,
   $\mu_t:= \L_{X_t}$ is a fixed point of $\Phi$ in $\C_V^\gg(\bar O)$. Since $\Phi$ has a unique fixed point, this implies the (weak) uniqueness
   of \eqref{E2}. Finally, by Gronwall's inequality, \eqref{KK2} follows from \eqref{P2}   for $X_t^\mu=X_t$ and $\mu_t:=\L_{X_t}$, where $\mu$ is the unique fixed point of $\Phi$. 
\end{proof}

  \paragraph{Acknowledgement.} The author would like to thank the referees for helpful comments and corrections.

\end{document}